\documentclass[a4paper, reqno]{amsart}          
\usepackage{amsmath, amssymb, latexsym, graphicx}

\newtheorem{thm}{Theorem}
\newtheorem{conj}{Conjecture}
\newtheorem{lem}{Lemma}
\newtheorem{cor}{Corollary}

\newcommand{\vol}[1]{\mathrm{vol}(#1)}
\newcommand{\R}{\mathbb{R}}
\newcommand{\norm}[1]{\left\|#1\right\|}
\newcommand{\vnorm}[1]{\|#1\|}
\newcommand{\bignorm}[1]{\bigl\|#1\bigr\|}
\newcommand{\abs}[1]{\left|#1\right|}

\begin{document}


\title{Cardinalities of $k$-distance sets in Minkowski spaces}
\date{1997}
\author{K.\ J.\ Swanepoel}
\address{Department of Mathematics and Applied Mathematics \\
University of Pretoria \\ Pretoria 0002 \\ South Africa \\ 
\texttt{konrad@math.up.ac.za}}


\maketitle
\begin{abstract}
A subset of a metric space is a {\em $k$-distance set\/} if there are exactly 
$k$ non-zero distances occuring between points. 
We conjecture that a $k$-distance set in a $d$-dimensional Banach space 
(or {\em Minkowski space}), contains at most $(k+1)^d$ points, with equality 
iff the unit ball is a parallelotope.
We solve this conjecture in the affirmative for all $2$-dimensional spaces and 
for spaces where the unit ball is a parallelotope.
For general spaces we find various weaker upper bounds for $k$-distance sets.
\end{abstract}


\section{Introduction}\label{Intro}
A subset $S$ of a metric space is a {\em $k$-distance set\/} if there are 
exactly $k$ non-zero distances occuring between points of $S$. 
We also call a $1$-distance set an {\em equilateral set.}
In this paper we find upper bounds for the cardinalities of $k$-distance sets 
in {\em Minkowski spaces}, i.e.\ finite-dimensional Banach spaces (see 
Theorems~\ref{thA} to \ref{Up}), and make a conjecture concerning tight upper 
bounds.

In Euclidean spaces $k$-distance sets have been studied extensively; see e.g.\ 
\cite{Erdos,Erdos2,Kelly,Golomb,BB,BBS,Blokhuis,Blokhuis2,Beck,Chung,CEGSW,CST,HP,Sz}, 
and the books \cite{Pach} and \cite[sections F1 and F3]{CFG}.

For general $d$-dimensional Minkowski spaces it is known that the maximum 
cardinality of an equilateral set is $2^d$, with equality iff the unit ball of 
the space is a parallelotope, and that if $d\geq 3$, there always exists an 
equilateral set of at least $4$ points \cite{Petty}.
It is unknown whether there always exists an equilateral set of $d+1$ points; 
see \cite{LM,Morgan} and \cite[p.\ 129, p.\ 308 problem 4.1.1]{Thompson}.
However, Brass \cite{Brass2} recently proved that for each $n$ there is a 
$d=d(n)$ such that any $d$-dimensional Minkowski space has an equilateral set 
of at least $n$ points.
See \cite{Guy} for problems on equilateral sets in $\ell_p$ spaces.

Equilateral sets in Minkowski spaces have been used in \cite{LM} to construct 
energy-minimizing cones over wire-frames.
See also \cite{Morgan}.

As far as we know, $k$-distance sets for $k\geq 2$ have not been studied in 
spaces other than euclidean.

Our main results are the following.

\begin{thm}  \label{thA}
If the unit ball of a $d$-dimensional Minkowski space is a parallelotope, then 
a $k$-distance set in $X$ has cardinality at most $(k+1)^d$.
This bound is tight.
\end{thm}

\begin{thm}\label{Cor1}
Given any set $S$ of $n$ points in a $d$-dimensional Minkowski space with a 
parallelotope as unit ball, there exists a point in $S$ from which there are 
at least $\lceil n^{1/d}\rceil-1$ distinct non-zero distances to points in $S$.
This bound is tight.
\end{thm}

\begin{thm}  \label{thB}
The cardinality of a $k$-distance set in a $2$-dimensional Minkowski space is 
at most $(k+1)^{2}$, with equality iff the space has a parallelogram as unit 
ball.
\end{thm}

\begin{thm}  \label{Cor2}
Given any set of $n$ points in a $2$-dimensional Minkowski space, there exists 
a point in $S$ from which there are at least $\lceil n^{1/2}\rceil-1$ distinct 
non-zero distances to points in $S$.
\end{thm}

\begin{thm}  \label{Up}
The cardinality of a $k$-distance set in a $d$-dimensional Minkowski space is 
at most $\min(2^{kd}, (k+1)^{(11^{d}-9^{d})/2})$.
\end{thm}

In the light of Theorems~\ref{thA} and \ref{thB} and the results of 
\cite{Petty}, we make the following

\begin{conj}
The cardinality of a $k$-distance set in any $d$-dimensional Minkowski space 
is at most $(k+1)^{d}$, with equality iff the unit ball is a parallelotope.
\end{conj}

As mentioned above, \cite{Petty} shows that this conjecture is true for $k=1$.
By Theorem~\ref{thB} the conjecture is true if $d=2$, and by Theorem~\ref{thA} 
if the unit ball is a parallelotope.

In the sequel, $(\R^d, \norm{\cdot})$ is a $d$-dimensional Minkowski space 
with norm $\norm{\cdot}$, $B(x,r)$ is the closed ball with centre $x$ and 
radius $r>0$, and $B:= B(0,1)$ the {\em unit ball}\/ of the space.
Recall that two $d$-dimensional Minkowski spaces are isometric iff their unit 
balls are affinely equivalent (by the Mazur-Ulam Theorem; see e.g.\ 
\cite[Theorem 3.1.2]{Thompson}).
In particular, a Minkowski space has a parallelotope as unit ball iff it is 
isometric to $(\R^d, \norm{\cdot}_{\infty})$, where 
$\norm{(\lambda_{1}, \lambda_{2},\dots,\lambda_{d})}_{\infty}:= 
\max_{i=1,\dots,d}\abs{\lambda_{i}}$.

We define a {\it cone} (or more precisely, an {\em acute cone}) $P$ to be a convex set in 
$\R^d$ that is positively homogeneous (i.e., for any $x\in P$ and 
$\lambda\geq 0$ we have $\lambda x\in P$) and satisfies $P\cap(-P)=\{0\}$.
Recall that such a cone defines a partial order on $\R^d$ by 
$x\leq y \iff y-x\in P$.

We denote the cardinality of a set $S$ by $\#S$.

For measurable $S\subseteq\R^{d}$, let $\vol{S}$ denote the Lebesgue measure 
of $S$.
For later reference we state Lyusternik's version of the Brunn-Minkowski 
inequality (see \cite[Theorem 8.1.1]{BZ}).

\begin{lem}
If $A,B\subseteq\R^{d}$ are compact, then
\[
  \vol{A+B}^{1/d} \geq \vol{A}^{1/d}+\vol{B}^{1/d}.
\]
If equality holds and $\vol{A}, \vol{B} > 0$, then $A$ and $B$ are convex 
bodies such that $A = v + \lambda B$ for some $\lambda >0$ and $v\in \R^d$.\qed
\end{lem}

\section{Proofs}

\begin{proof}[Proof of Theorem~\ref{thA}]
We may assume without loss of generality that the space is 
$(\R^d,\norm{\cdot}_\infty)$.
We introduce partial orders on $\R^d$ following Blokhuis and Wilbrink 
\cite{BlW}.
For each $i=1,\dots,d$, let $\leq_{i}$ be the partial order with cone
\[
  P_{i}=\bigl\{(\lambda_{1},\dots,\lambda_{d})\in\R^d: 
               \max_{j=1,\dots,d}\abs{\lambda_j}=\lambda_{i}\bigr\}.
\]

For each $x$ in a $k$-distance set $S$, let $h_{i}(x)$ be the length of the 
longest descending 
$\leq_{i}$-chain starting with $x$, i.e.\ $h_{i}(x)$ is the largest $h$ such 
that there exist $x_{1},x_{2},\dots,x_{h}\in S$ for which 
$x >_{i} x_{1} >_{i} x_2 >_i \dots >_{i} x_{h}$.

Since $\bigcup_{i=1}^{d}(P_{i}\cup -P_{i}) = \R^d$, for all 
distinct $x,y\in\ell_{\infty}^{d}$ there exists $i$ such that $x <_{i} y$ or 
$y <_{i} x$.
Exactly as in \cite{BlW}, it follows that the mapping 
$x\mapsto (h_{1}(x),\dots,h_{d}(x))$ is injective, and thus 
$\#S\leq (h+1)^{d}$, where
\[
  h:=\max_{x\in S, i=1,\dots, d} h_{i}(x).
\]

It remains to show that $h\leq k$.
Suppose not.
Then for some $x\in S$ and some $i$ there exist $x_{1},\dots,x_{k+1}\in S$ 
such that $x >_{i} x_{1} >_{i} \dots >_{i} x_{k+1}$.
Since $S$ is a $k$-distance set, 
$\norm{x-x_{m}}_{\infty}=\norm{x-x_{n}}_{\infty}$ for some 
$1\leq m < n \leq k+1$.
Also, $x-x_{m}, x-x_{n}\in P_{i}$.
Now note that if $\norm{a}_{\infty}=\norm{b}_{\infty}$ with 
$a,b\in P_{i}, a\neq b$, then $a$ and $b$ are $\leq_{i}$-incomparable; 
in particular, $b-a\not\in P_{i}$.
Therefore, $x_{m}-x_{n}\not\in P_{i}$, a contradiction.

The set $\{0,1,\dots,k\}^d$ is a $k$-distance set of cardinality $(k+1)^d$.
Note that it is not difficult to see that in fact the only $k$-distance sets 
of cardinality $(k+1)^d$ are of the form $S=a+\lambda\{0,1,\dots,k\}^d$ for 
some $a\in\R^d$ and $\lambda > 0$.
\qed
\end{proof}

\begin{proof}[Proof of Theorem~\ref{Cor1}]
Consider the mapping $x\mapsto (h_1(x),\dots,h_d(x))$ in the proof of 
Theorem~\ref{thA}.
If $h$ is the length of the longest $\leq_i$-chain over all $i$, then 
$n\leq (h+1)^d$.
Thus there is a $\leq_i$-chain $x_0>_i x_1>_i \dots >_i x_h$ of length 
$h\geq \lceil n^{1/d}\rceil-1$.
By the last paragraph of the proof of Theorem~\ref{thA}, the distances 
$\rho(x_0,x_j), j=1,\dots,h$ are all distinct.

Any $S\subseteq\R^d$ such that
\[
\{0,1,\dots,\lceil n^{1/d}\rceil-2\}^d\subsetneq S 
\subseteq \{0,1,\dots,\lceil n^{1/d}\rceil-1\}^d,
\]
has exactly $\lceil n^{1/d}\rceil-1$ distinct distances in the norm 
$\norm{\cdot}_\infty$.
\qed
\end{proof}

The following corollary is easily gleaned from the proof of Theorem~\ref{thA}.

\begin{cor}  \label{Cor}
Suppose that $\{P_{i}: i\in I\}$ is a family of cones in a Minkowski space 
$(\R^{d},\norm{\cdot})$ satisfying
\begin{equation}  \label{one}
  \bigcup_{i\in I} (P_{i}\cup - P_{i}) = \R^{d},
\end{equation}
and
\begin{equation}  \label{two}
  \forall\, i\in I \: \forall\, \text{distinct } x,y\in P_{i},
    \text{ if } \norm{x}=\norm{y} \text{ then } \pm(x-y)\not\in P_{i}.
\end{equation}
Then a $k$-distance set in $(\R^{d},\norm{\cdot})$ has cardinality at most 
$(k+1)^{\#I}$.
\qed
\end{cor}

\begin{lem}\label{metriclemma}
Let $S$ be a $k$-distance set in a metric space $(X,\rho)$ with distances 
$\rho_{1} < \rho_{2} < \dots < \rho_{k}$.
If $\rho_{k}/\rho_{1} > 2^{k-1}$, then for some $i=1, \dots, k-1$, the relation
\[
  x\sim_{i} y \iff \rho(x,y) \leq \rho_{i}
\]
is an equivalence relation.
\end{lem}
\begin{proof}
The relation $\sim_{i}$ is reflexive and symmetric.
If it is not transitive, there exist $x,y,z\in S$ such that 
$\rho(x,y),\rho(y,z) \leq \rho_{i}$ and $\rho(x,z) > \rho_{i}$.
Thus $\rho_{i+1}\leq \rho(x,z) \leq \rho(x,y)+\rho(y,z) \leq 2\rho_{i}$.
If this holds for all $i=1, \dots, k-1$, we obtain 
$\rho_{k}\leq 2^{k-1}\rho_{1}$.
\qed
\end{proof}

\begin{lem}  \label{Up1}
The cardinality of a $k$-distance set in a $d$-dimensional Minkowski space is 
at most $2^{kd}$.
\end{lem}
\begin{proof}
Let $\{x_{1}, \dots, x_{m}\}$ be a $k$-distance set with distances 
$\rho_{1} < \rho_{2} < \dots < \rho_{k}$.
Set $V := \bigcup_{i=1}^{m} B(x_{i}, \rho_{1}/2)$.
Then we have 
\begin{equation}  \label{voleq}
\vol{V} = m(\rho_{1}/2)^{d}\vol{B}.
\end{equation}
Also, $V-V\subseteq B(0, \rho_{k}+\rho_{1})$, since if $x,y\in V$, there exist 
$i$ and $j$ such that $\norm{x-x_{i}}\leq \rho_{1}/2$, 
$\norm{y-x_{j}}\leq \rho_{1}/2$.
Thus
\[
  \norm{x-y} \leq \norm{x-x_{i}} + \norm{x_{i} - x_{j}} + \norm{x_{i}-x_{j}} 
  \leq \rho_{1} + \rho_{k}.
\]
Therefore, 
\begin{equation}
\label{voleq2}
\vol{V-V}\leq (\rho_{1}+\rho_{k})^{d}\vol{B}.
\end{equation}
Substituting \eqref{voleq} and \eqref{voleq2} into the Brunn-Minkowski 
inequality
\begin{equation}  \label{bmineq}
  \vol{V-V}^{1/d}\geq\vol{V}^{1/d}+\vol{-V}^{1/d},
\end{equation}
we obtain
$\rho_{1}+\rho_{k}\geq m^{1/d}\rho_{1}$, and $m\leq (1+\rho_{k}/\rho_{1})^{d}$.

If $1+\rho_{k}/\rho_{1}\leq 2^k$, there is nothing to prove.
Otherwise, $\rho_{k}/\rho_{1} > 2^k-1 \geq 2^{k-1}$, and by 
Lemma~\ref{metriclemma}, $x\sim_{i} y \iff \rho(x,y) \leq \rho_{i}$ is an 
equivalence relation for some $i=1,\dots,k-1$.
By induction on $k$ we obtain that each equivalence class, being an 
$i$-distance set, has at most $2^{id}$ points.
By choosing a representative from each equivalence class, we obtain a 
$(k-i)$-distance set with at most $2^{(k-i)d}$ points.
Therefore, $m\leq 2^{id}2^{(k-i)d} = 2^{kd}$.
\qed
\end{proof}

In the proof of Theorem~\ref{thB}, we need the following geometric lemma, 
which is a modification of \cite[corollary 3.2.6]{Thompson} in $2$ dimensions.

\begin{lem}  \label{Auerbach}
Let $B_1$ be the convex hull of $\{(\pm 1,0), (0,\pm 1)\}$ and $B_\infty$ the 
square $[-1,1]^2$.
For any symmetric convex disc $C$ in $\R^{2}$ there exists an invertible linear
transformation taking $C$ to $C'$ such that 
$B_1\subseteq C'\subseteq B_\infty$ and such that any straight-line segment 
contained in the boundary of $C'$ lies completely in one of the four 
coordinate quadrants.
\end{lem}

\begin{proof}
We consider all triangles with vertices $0,x,y$, where $x$ and $y$ are on the 
boundary of $C$.
By compactness there exist $x_{0}$ and $y_{0}$ such that the area of the 
triangle is a maximum.
Then $\{x_{0}+\lambda y_{0}: \lambda\in\R\}$ is a support line of $C$ at 
$x_{0}$, since otherwise we can replace $x_{0}$ by a point on the side of the 
line opposite $0$ to enlarge the area of the triangle.
Similarly, $\{y_{0}+\lambda x_{0}: \lambda\in\R\}$ is a support line of $C$ at 
$y_{0}$.
Since $C$ is symmetric, it follows that $C$ is contained in the parallelogram 
$\{\lambda x_{0}+\mu y_{0}:-1\leq \lambda,\mu\leq 1\}$.
See Figure~\ref{auerfig}.
\begin{figure}
  \begin{center}
    \includegraphics{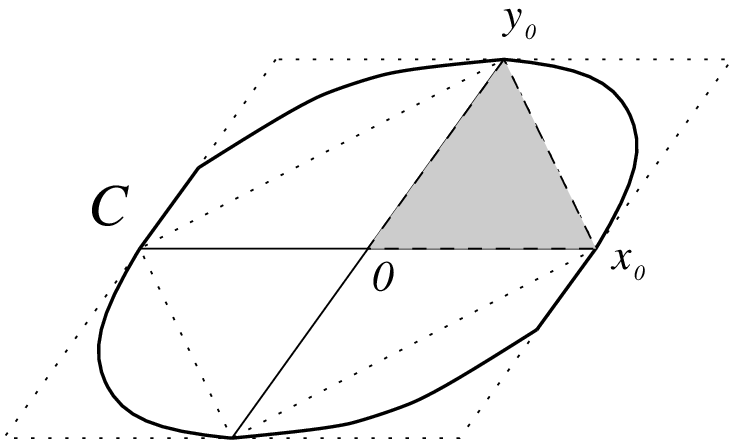}
  \end{center}
\caption{}\label{auerfig}
\end{figure}

If $x_{0}$ is an interior point of a straight-line segment contained in the 
boundary of $C$, we may shift $x_{0}$ to a boundary point of such a segment, 
without changing the area of the triangle.
Thus $C$ is still contained in a parallelogram as above.
A similar remark holds for $y_{0}$.
We now apply the linear transformation sending $x_{0}$ and $y_{0}$ to the 
standard unit vectors $e_1$ and $e_2$, respectively (see 
Figure~\ref{auerfig2}).
\qed
\end{proof}

\begin{proof}[Proof of Theorem~\ref{thB}]
We have to find two cones $P_{1}$ and $P_{2}$ satisfying \eqref{one} and 
\eqref{two} of Corollary~\ref{Cor}.
By Lemma~\ref{Auerbach} we may replace the space by an isometric space 
$(\R^{2},\norm{\cdot})$ such that the unit ball $B$ of $\norm{\cdot}$ lies 
between $B_1$ and $B_\infty$, and such that any straight-line segment 
contained in the boundary of the unit ball lies completely in a quadrant of 
the plane.

We provisionally let $P_{1}$ be the closed first quadrant, and $P_{2}$ the 
closed second quadrant.
See Figure~\ref{auerfig2}.
\begin{figure}
  \begin{center}
    \includegraphics{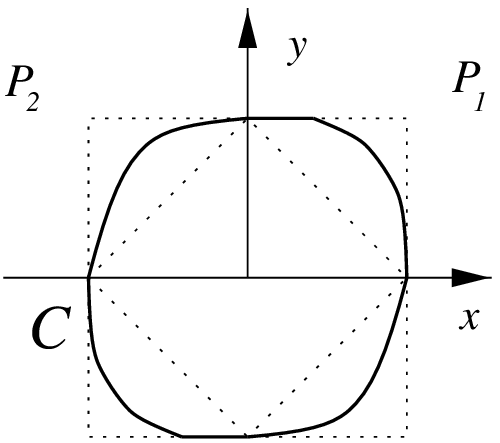}
  \end{center}
\caption{}\label{auerfig2}
\end{figure}
Then \eqref{one} is satisfied.
The only way that \eqref{two} could fail is if there is a straight-line 
segment contained in the boundary of the unit ball parallel to either the 
x-axis or the y-axis, lying in $P_{1}$ or $P_{2}$.
If there is a segment in the boundary of the unit ball in $P_{1}$ parallel to 
the x-axis, say, we remove the positive x-axis $\{(\lambda,0): \lambda > 0\}$ 
from $P_{1}$.
If in this case there were another straight-line segment in the boundary 
parallel to the x-axis in $P_{2}$, then there would be a straight-line segment 
in the boundary lying in the first and second quadrants, giving a 
contradiction.
Thus we do not have to remove the negative x-axis from $P_{2}$, and \eqref{one}
is still satisfied.
We do the same thing for segments parallel to the y-axis, and for $P_{2}$.
In the end, the modified $P_{1}$ and $P_{2}$ satisfy \eqref{one} and 
\eqref{two}, and we deduce $\#S\leq (k+1)^{2}$ from Corollary~\ref{Cor}.

If equality holds, then the mapping $x\mapsto (h_{1}(x),h_{2}(x))$
in the proof of Theorem~\ref{thA} is a bijection from $S$ to 
$\{0,\dots,k\}^{2}$.
We now denote a point $x\in S$ by $p_{i,j}$, where 
$(i,j) = (h_{1}(x),h_{2}(x))$.

Suppose that two of the distances $\norm{p_{0,i}-p_{0,0}}$ ($i=1,\dots,k$) are 
equal, say $\norm{p_{0,i}-p_{0,0}} = \norm{p_{0,j}-p_{0,0}}$ with $0<i<j$.
Then, since $p_{0,j} >_2 p_{0,i} >_2 p_{0,0}$, we have 
$p_{0,i}-p_{0,0}, p_{0,j}-p_{0,0} \in P_2$, contradicting \eqref{two}.

It follows that the distances $\norm{p_{0,i}-p_{0,0}}$, $i=1,\dots,k$ are 
distinct, and thus are exactly the $k$ different distances in increasing order.
Similarly, the distances $\norm{p_{0,i}-p_{0,1}}$, $i=2,\dots,k$ are in 
increasing order.
If $\norm{p_{0,k}-p_{0,1}}=\rho_{k}$, the three points 
$p_{0,0},p_{0,1},p_{0,k}$ again contradict \eqref{two}.
Thus these distances are $\rho_{1},\dots, \rho_{k-1}$ in increasing order, etc.
In the end we find that $\norm{p_{0,i+1}-p_{0,i}}=\rho_{1}$ for all $i$.
Thus $\rho_{k}\leq k\rho_{1}$, by the triangle inequality.
Using the Brunn-Minkowski inequality as in the proof of Lemma~\ref{Up1}, we 
find that equality holds in \eqref{bmineq} and \eqref{voleq2}, implying that 
for $V := \bigcup_{i=1}^{\#S} B(x_{i}, \rho_{1}/2)$ we have 
$V-V = B(0,\rho_{k}+\rho_{1})$, and $V-V$ and $V$ are homothetic. 
Thus $V$ is a ball that is perfectly packed by smaller balls.
By a result of \cite{Groemer}, this implies that the unit ball is a 
parallelogram.
\qed
\end{proof}

\begin{proof}[Proof of Theorem~\ref{Cor2}]
Follows from the proof of Theorem~\ref{thB} in the same
way that Theorem~\ref{Cor1} follows from Theorem~\ref{thA}.
\qed
\end{proof}

\begin{proof}[Proof of Theorem~\ref{Up}]
Lemma~\ref{Up1} already gives part of the theorem.
For the remaining part we apply Corollary~\ref{Cor}.
In order for a cone $P$ to satisfy \eqref{two}, it is sufficient that 
\begin{equation}  \label{condition}
\forall\, a,b\in P: \text{ if } \norm{a}=\norm{b}=1, 
\text{ then } \norm{a-b}<1.
\end{equation}
To see this, suppose that $P$ does not satisfy the condition in \eqref{two}, 
i.e.\ there exist distinct $x,y\in P$ such that $\norm{x}=\norm{y}$ and 
$y-x\in P$.
Let $a:= \norm{x}^{-1}x$, $b:= \norm{y}^{-1}y$, $c:= \norm{y-x}^{-1}(y-x)$, 
and $0<\lambda := \norm{x}/(\norm{y-x}+\norm{x}) < 1$.
Then $a=(1-\lambda)(a-c)+\lambda b$, and 
\[
  1=\norm{a}\leq(1-\lambda)\norm{a-c}+\lambda\norm{b}
   =(1-\lambda)\norm{a-c}+\lambda,
\]
implying $\norm{a-c}\geq 1$.

In order for \eqref{one} to be satisfied too, we need a cover of the unit 
sphere by sets such that, if they are extended to positive cones, are convex.

We do this with the following construction:
Let $C=\{c_{1},c_{2},\dots,c_{m}\}$ be a maximal set of unit vectors satisfying
$\norm{c_{i} - c_{j}}, \norm{c_{i} + c_{j}} \geq\tfrac{1}{5}$ for all 
$1\leq i<j\leq m$.
Then for any unit vector $x$ there exists $i$ such that 
$\norm{c_{i} - x} < \tfrac{1}{5}$ or $\norm{c_{i} + x} < \tfrac{1}{5}$.
For $i=1,\dots,m$, let $P_{i}$ be the cone generated by
\[
  Q_{i}:=\bigl\{x\in\R^{d}:\norm{x}=1,\norm{c_{i}-x}<\tfrac{1}{5}\bigr\},
\]
i.e.\ $P_{i}:=\{\sum_{j}\lambda_{j}x_{j}:\lambda_{j}\geq 0, x_{j}\in Q_{i}\}$.
Then the $P_{i}$'s satisfy \eqref{one} by the maximality of $C$.
Each $P_{i}$ satisfies \eqref{condition}:
Let $\sum_{j}\lambda_{j}x_{j}\in P_{i}$, where 
$\lambda_{j}\geq 0, \|x_{j}\|=1, \norm{c_{i}-x_{j}}<\tfrac{1}{5}$ and 
$\norm{\sum\lambda_{j}x_{j}}=1$.
Then 
\begin{align*}
  \bignorm{c_{i}-\sum_{j}\lambda_{j}x_{j}} 
  & = \bignorm{\sum_{j}\lambda_{j}(c_{i}-x_{j}) 
       + (1-\sum_{j}\lambda_{j})c_{i}} \\
  & < \sum_{j}\lambda_{j}/5 -1 + \sum_{j}\lambda_{j} 
           \quad\text{(since } \sum_{j}\lambda_{j}\geq 1) \\
  & = \tfrac{6}{5}\sum_{j}\lambda_{j} - 1.
\end{align*}
Also, since
\[
  1+\sum_{j}\lambda_{j}/5 
  > \bignorm{\sum_{j}\lambda_{j}x_{j}} + \sum_{j}\norm{\lambda_{j}x_{j}
                -\lambda_{j}c} 
  \geq \sum_{j}\lambda_{j}\norm{c}=\sum_{j}\lambda_{j},
\]
we obtain $\sum_{j}\lambda_{j} < \tfrac{5}{4}$, and 
$\vnorm{c_{i}-\sum_{j}\lambda_{j}x_{j}} 
< \tfrac{6}{5}\cdot\frac{5}{4} - 1 
= \tfrac{1}{2}$.

A volume argument gives the upper bound for $\# C$:
The balls
\[
  B(0,\tfrac{9}{10}), B(\pm c_{i},\tfrac{1}{10}), i=1,\dots m
\]
have disjoint interiors and are contained in the ball $B(0,\tfrac{11}{10})$.
Therefore,
\[
  (\tfrac{9}{10})^{d}\vol{B} + 2m(\tfrac{1}{10})^{d}\vol{B} 
  \leq (\tfrac{11}{10})^{d}\vol{B},
\]
giving $m\leq \tfrac{1}{2} (11^{d}-9^{d})$.
\qed
\end{proof}

\section*{Acknowledgement}
This paper is part of the author's PhD thesis written under supervision of 
Prof.\ W. L. Fouch\'e at the University of Pretoria.
I thank the referees as well as Graham Brightwell for their suggestions on the 
layout of the paper.

\providecommand{\bysame}{\leavevmode\hbox to3em{\hrulefill}\thinspace}

\end{document}